\documentclass[11pt]{amsart}
\usepackage{amssymb,amscd,epsf,amsmath}

\begin{document}

\newtheorem{defi}{Definition}
\newtheorem{thm}{Theorem}
\newtheorem{prop}{Proposition}
\newtheorem{lemma}{Lemma}
\newtheorem{coro}{Corollary}
\newtheorem{quest}{Question}

\theoremstyle{remark}
\newtheorem*{remark}{Remark}
\newtheorem{numrem}{Remark}

\title[An Extremum Property of the Cross-Polytope]{An Extremum
Property Characterizing the\\
$n$-dimensional Regular Cross-Polytope}
\author{Wlodzimierz Kuperberg}
\date{\today}
\address{Department of Mathematics, Auburn University, Auburn, AL 36849-5310}
\email{kuperwl@auburn.edu}

\begin{abstract} In the spirit of the \textit{Genetics of the Regular Figures},
by L.~Fejes T\'oth \cite[Part~2]{lft}, we prove the following theorem: If $2n$
points are selected in the $n$-dimensional Euclidean ball $B^n$ so that the
smallest distance between any two of them is as large as possible, then the
points are the vertices of an inscribed regular cross-polytope. This
generalizes a result of R.A.~Rankin \cite{rankin} for $2n$ points on the
\emph{surface} of the ball. We also generalize, in the same manner, a theorem
of Davenport and Haj\'{o}s \cite{lapok1} on a set of $n+2$ points. As a
corollary, we obtain a solution to the problem of packing $k$ unit
$n$-dimensional balls $(n+2\le k\le 2n)$ into a spherical container of minimum
radius. \end{abstract}

\maketitle

\section{Introduction}

The {\it regular cross-polytope} is the dual to the $n$-dimensional Euclidean
cube. More directly, the regular cross-polytope can be described as the convex
hull of the union of $n$ mutually perpendicular line segments of equal length,
intersecting at the midpoint of each of them. Obviously, the regular
cross-polytope generated by perpendicular segments of lenth $d$ is inscribed in
a sphere of radius $r=d/2$, and the edge-length of the cross polytope is
$\sqrt2\, r$. Each of the $2^n$ facets of the cross-polytope is a regular
$(n-1)$-dimensional simplex. The only distances between vertices of the regular
cross-polytope inscribed in a ball of radius $r$ are $\sqrt2\, r$ and $2r$.

The main goal of this note is to prove the following metric characterization of
the regular cross-polytope in terms of an extremum property of its vertex set: 

\begin{thm}\label{main} If $\,V$ is a $2n$-point subset $(n\ge 2)$ of the unit
ball in $n$-dimensional Euclidean space such that the shortest
distance between points in $V$ is as large as possible, then $V$ is the
set of vertices of a regular cross-polytope inscribed in the ball.
\end{thm}

This characterization of the $n$-dimensional regular cross-polytope by the
extremum property of the distances between its vertices fits well into L.~Fejes
T\'{o}th's theory of {\it genetics of regular figures}  \cite[Part~2]{lft}. The
theory is supported by a collection of examples illustrating how `` [ ...
]~\emph{regular arrangements are generated from unarranged chaotic sets by
ordering effect of an economy principle, in the widest sense of the word}\/''
\cite[Preface, p.~x]{lft}.

Theorem~\ref{main} is a generalization of a result obtained by R.A.~Rankin
\cite{rankin} in 1955. In Rankin's version it is assumed that the points lie on
the surface of the ball. Rankin states the result in a comment at the end of
the proof of one of his theorems \cite[p.~142]{rankin}. Before Rankin,
K.~Sch\"{u}tte and B.L.~van~der~Waerden \cite{sw} solved the $3$-dimensional
case, also with the assumption that the points lie on the ball's surface.

In a similar way we generalize a theorem stated by Davenport and Haj\'{o}s
\cite{lapok1} (proved in \cite{lapok2} and in~\cite{rankin}), concerning a set
of $n+2$ points. Again, in our generalization the restriction that the points
lie on the surface of the ball is removed:

\begin{thm}\label{modif} If $n+2$ points lie in the $n$-dimensional Euclidean
unit ball, then at least one of the distances between the points is smaller
than or equal to $\sqrt 2$.
\end{thm}

As a direct application of the above two theorems we obtain the following result
on packing balls in a spherical container:

\begin{thm}\label{pack} Let $S$ be a spherical container of minimum radius that
can hold $n+2$ nonoverlapping $n$-dimensional balls $(n\ge 2)$ of radius $1$
each. Then the radius of $S$ is $1+\sqrt2$, and there is enough room in $S$ to
hold as many as $2n$ such balls. Moreover, the packing configuration of $2n$
balls in $S$ is unique up to isometry, the balls' centers forming the set of
vertices of an $n$-dimensional cross-polytope. Hence $S$ cannot accommodate
$2n+1$ unit balls.
\end{thm}

\begin{remark} The problem of packing $k$ unit $n$-dimensional balls in a
spherical container of minimum radius $r=r(k,n)$ for $k\le n+1$ easily
reduces to the problem of distributing $k$ points on the \emph{surface} of the
ball so that the shortest distance between the points is as large as possible.
This problem has been solved by Rankin \cite{rankin}: the points are vertices
of a regular $(k-1)$-dimensional simplex inscribed in the ball and concentric
with it. Thus,
$$r(k,n)=1+\sqrt{2-\frac{2}{k}}\ \ {\rm for}\ \ k\le n+1,$$ 
while
according to Theorem \ref{pack}
$$r(k,n)=1+\sqrt2\ \ {\rm for}\ \ n+2\le k\le 2n,\ \ {\rm and}\ \
r(2n+1,n)>1+\sqrt2.$$
\end{remark}

\section{Notation and Preliminary Statements}

The cardinality of a set $A$ is denoted by ${\rm card}A$. $\mathbb R^n$ denotes
Euclidean (Cartesian) $n$-dimensional space with the usual inner product and
the metric thereby generated. The origin $(0,0,\ldots ,0)$ is denoted by $o$.
$B^n$ denotes the unit ball in $\mathbb R^n$, and $S^{n-1}$ denotes the unit
sphere, {\it i.e.,} the boundary of $B^n$. For $x\in S^{n-1}$, the hyperplane
containing the origin and normal to $x$ is denoted by $H_0(x)$, and $H(x)$
denotes the half-space of $H_0(x)$ opposite to $x$. In other words,
$$H_0(x)=\left\{y\in \mathbb R^n : xy=0\right\},$$
and
$$H(x)=\left\{y\in \mathbb R^n : xy\le 0\right\}.$$

Also, with each $x\in S^{n-1}$ we associate the set
$$C(x)=\left\{y\in B^n : {\rm dist}(x,y)\ge\sqrt2\right\},$$
called \emph{the crescent determined by} $x$. (For $n=2$, the shape of $C(x)$
resembles a crescent.)

Further, for $A\subset \mathbb R^n$, ${\rm Conv}A$ denotes the convex hull of
$A$ and ${\rm Lin}A$ denotes the linear hull of $A$, that is, the smallest
linear subset of $\mathbb R^n$ containing $A$. The interior of $A$ is denoted by ${\rm
Int}A$ and ${\rm Int}_LA$ denotes the interior of $A$ relative to ${\rm Lin}A$.

Next, we state a few propositions, simple enough to have their proofs omitted.

\begin{prop} $C(x)\subset H(x)$ for every $x\in S^{n-1}$.
\end{prop}

\begin{prop} $C(x)\cap H_0(x)$ is a great $(n-2)$-sphere of $S^{n-1}$.
Specifically, $C(x)\cap H_0(x)=H_0(x)\cap S^{n-1}$.
\end{prop}

\begin{prop} If $x\in S^{n-1}$, then for every point $y$ in the closed half-ball
$B^n\setminus {\rm Int}H(x)$, we have ${\rm dist}(x,y)\le\sqrt2$.
\end{prop}

\begin{prop} Suppose $A\subset S^{n-1}$. If\ \ $o\in{\rm
Int}({\rm Conv}A)$, then\ \ $\displaystyle\bigcap_{x\in A}H(x)=\{ o\}$.
\end{prop}

The above proposition is immediately generalized to:

\begin{prop} Suppose $A\subset S^{n-1}$ and $k=\dim {\rm Lin}A$. If\
\ $o\in{\rm Int}_L({\rm Conv}A)$, then\ \ $\displaystyle\bigcap_{x\in A}H(x)$\ \ is
the $(n-k)$-dimensional linear subspace of $\mathbb R^n$ normal to ${\rm Lin}A$. In
particular,\ \ $\displaystyle\bigcap_{x\in A}H(x)=\bigcap_{x\in A}H_0(x)$.
\end{prop}

The above propositions imply directly the following intersection properties of
the crescents determined by a subset of $S^{n-1}$:

\begin{prop} Suppose $A\subset S^{n-1}$. If\ \ $o\in{\rm
IntConv}A$, then\ \ $\displaystyle\bigcap_{x\in A}C(x)=\emptyset$.
\end{prop}

\begin{prop} Suppose $A\subset S^{n-1}$ and $k=\dim{\rm Lin}A$. If\ \
$o\in{\rm Int}_L({\rm Conv}A)$, then\ \ $\displaystyle\bigcap_{x\in A}C(x)$\ \
is the $(n-k-1)$-dimensional great sphere of $S^{n-1}$ lying in the linear subspace
normal to ${\rm Lin}A$.
\end{prop}

Every bounded subset $A$ of $\mathbb R^n$ containing at leat two points is
contained in a (unique) ball of minimum radius. The radius of the smallest ball
containing $A$ is called the {\it circumradius} of $A$ and is denoted by $r(A)$.
If $A$ is compact, then $A$ contains a {\it finite} subset of the same
circumradius as $A$. Obviously, among such finite subsets of $A$ there is one
(not necessarily unique) of minimum cardinality.

\begin{prop} Suppose $F$ is a finite subset of $B^n$ with $r(F)=1$. If $r(F')<1$
for every proper subset $F'$ of $F$, then $\dim{\rm Lin}F={\rm
card}F-1$, $F$ lies on $S^{n-1}$ and $o\in{\rm Int}_L({\rm Conv}F)$. 
\end{prop}

\section{Proofs}

\begin{proof}[Proof of Theorem \ref{modif}] Let $P=\{p_1,
p_2,\ldots,p_{n+2}\}$ be a subset of $B^n$. Since our goal is to show that one of
the distances between points in $P$ is smaller than or equal to
$\sqrt2$, we may assume that the circumradius of $P$ is 1, expanding $P$
homothetically if needed. By this assumption, it follows that $o\in{\rm
Conv}P$. By a well-known theorem of Carath\'{e}odory, $o$ lies in the convex
hull of an
$(n+1)$-element subset of $P$, say $o\in{\rm Conv}P_1$, where
$P_1=P\setminus\{p_1\}$. This implies that $P_1$ is not contained in ${\rm
Int}H(p_1)$, {\it i.e.}, some point of
$P_1$, say $p_2$, lies in the half-space complementary to $H(p_1)$. Therefore
$p_2$ lies in the closed half-ball $B^n\setminus{\rm Int}H(p_1)$, which
implies ${\rm dist}(p_1,p_2)\le\sqrt2$. \end{proof}

\begin{proof}[Proof of Theorem \ref{main}] As we noted before, the shortest
distance between the $2n$ vertices of the regular cross-polytope inscribed in
$B^n$ is $\sqrt2$. It follows that the distance between any two points in $V$
is greater than or equal to $\sqrt2$. Observe that $r(V)=1$, for otherwise all
distances between points of $V$ could be enlarged by expanding $V$
homothetically. Let $V_0\subset V$ be a set of minimum cardinality among all
subsets of $V$ whose circumradius is~$1$, and let $k={\rm card}V_0$. Obviously,
$k\ge 2$, and by the theorem of Carath\'{e}odory, $k\le n+1$. We now proceed
inductively:

$1^\circ$ If $n=2$, then $V$ is a four-point subset of the unit disk $B^2$, and
it could well be left to the reader to show that $V$ is the set of vertices
of an inscribed square. Nevertheless, we present here the following argument
since it will also serve as an illustration to our inductive step $2^\circ$.
By Proposition 8, $V_0$ is a subset of $S^1$, $\dim{\rm Lin}V_0=k-1$ is
either $2$ or $1$ and $o\in {\rm Int}_L({\rm Conv}V_0)$. Now, every point of
$V\setminus V_0$ lies in\ \ $\displaystyle\bigcap_{x\in V_0}C(x)$. Hence, by
Proposition~7, the set $V\setminus V_0$ is contained in a ``great
$m$-dimensional sphere'' of $S^1$, where $m=2-k$. Since the set $V\setminus V_0$
is nonempty, $m$ cannot be negative. Thus $k=2$, and since $o\in {\rm
Int}_L({\rm Conv}V_0)$, the two points of $V_0$ are antipodes. The remaining
two pints of $V$ lie in the ``great $0$-dimensional sphere'' of $S^1$ determined
by the line normal to $L(V_0)$, which concludes this portion of the proof.

$2^\circ$ Assume that $n\ge3$ and that the conclusion of Theorem \ref{main} is
true for balls of dimension smaller than $n$. By Proposition 8, $V_0$ is a
subset of $S^{n-1}$, $\dim{\rm Lin}V_0=k-1$, and $o\in {\rm Int}_L({\rm
Conv}V_0)$. Now, every point of $V\setminus V_0$ lies in\ \
$\displaystyle\bigcap_{x\in V_0}C(x)$. Hence, by Proposition~7, the set
$V\setminus V_0$ is contained in a great $m$-dimensional sphere of $S^{n-1}$,
where $m=n-k$. By the inductive assumption, the cardinality of $V\setminus V_0$
cannot exceed $2(m+1)$, the number of vertices of the $(m+1)$-dimensional
cross-polytope. That means
$2n-k\le2(n-k+1)$, which implies $k\le2$. Therefore $V_0$ consists of a pair of
antipodes of $S^{n-1}$, and the remaining $2n-2$ points of $V$ lie on the great
$(n-2)$-dimensional sphere of $S^{n-1}$, in a hyperplane perpendicular to the
line containing $V_0$. By the inductive assumption, $V\setminus V_0$ is the set
of vertices of an $(n-1)$-dimensional cross-polytope, hence $V$ is the set of
vertices of an $n$-dimensional cross-polytope inscribed in $S^{n-1}$.
\end{proof}


\begin{thebibliography}{[AB]}

\bibitem[SW]{sw} K.~Sch\"{u}tte and B.L.~van~der~Waerden,  \emph{
Auf welcher Kugel haben 5, 6, 7, 8 oder 9 Punkte mit Mindestabstand 1 Platz?\/}, 
Math.~Ann {\bf 123} (1951),  139-144.

\bibitem[DH]{lapok1} H.~Davenport and G.~Haj\'{o}s,  \emph{
Problem 35\/} (in Hungarian),  Mat. Lapok {\bf 2} (1951),  68.

\bibitem[AS]{lapok2} J.~Acz\'{e}l and T.~Szele,  \emph{Solutions to
Problem 35\/} (in Hungarian),  Mat. Lapok {\bf 3} (1952),  94-95.

\bibitem[R]{rankin} R.A.~Rankin,  \emph{
The closest packing of spherical caps in $n$ dimensions\/},  Proc.~Glasgow
Math.~Assoc.~{\bf 2} (1955),  139-144.

\bibitem[FT]{lft} L.~Fejes T\'{o}th, \emph{Regular Figures\/}, Pergamon Press,
1964.

\end{thebibliography}
\end{document}